\documentclass[12pt]{article}

\usepackage{amsmath,amsthm}
\usepackage{amssymb}
\usepackage{color}
\usepackage{xspace}
\usepackage[colorlinks=true,
linkcolor=green,
filecolor=brown,
citecolor=blue]{hyperref}

\def\gl{ground level\xspace}

\def\v{\vert}
\def\a{\ensuremath{\mathcal A}\xspace}
\def\bb{\ensuremath{\mathcal B}\xspace}
\def\f{\ensuremath{\mathcal F}\xspace}

\def\mbf#1{\mathchoice{\hbox{\boldmath $\displaystyle #1$}}
        {\hbox{\boldmath $\textstyle #1$}}
        {\hbox{\boldmath $\scriptstyle #1$}}
        {\hbox{\boldmath $\scriptscriptstyle #1$}}} 

\catcode`\@=11

\thicklines
\newskip\Einheit \Einheit=.6cm
\newcount\xcoord \newcount\ycoord
\newdimen\xdim \newdimen\ydim \newdimen\PfadD@cke \newdimen\Pfadd@cke
\def\PfadDicke#1{\PfadD@cke#1 \divide\PfadD@cke by2 
\Pfadd@cke\PfadD@cke \multiply\PfadD@cke by2}
\long\def\LOOP#1\REPEAT{\def\BODY{#1}\ITERATE}
\def\ITERATE{\BODY \let\next\ITERATE \else\let\next\relax\fi \next}
\let\REPEAT=\fi
\def\Punkt{\hbox{\raise-2pt\hbox to0pt{\hss\scriptsize$\bullet$\hss}}}

\def\DuennPunkt(#1,#2){\unskip
  \raise#2 \Einheit\hbox to0pt{\hskip#1 \Einheit
          \raise-1.5pt\hbox to0pt{\hss\tiny$\bullet$\hss}\hss}}
		  
\def\NormalPunkt(#1,#2){\unskip
  \raise#2 \Einheit\hbox to0pt{\hskip#1 \Einheit
          \raise-3pt\hbox to0pt{\hss\large$\bullet$\hss}\hss}}
\def\DickPunkt(#1,#2){\unskip
  \raise#2 \Einheit\hbox to0pt{\hskip#1 \Einheit
          \raise-4pt\hbox to0pt{\hss\Large$\bullet$\hss}\hss}}
\def\Kreis(#1,#2){\unskip
  \raise#2 \Einheit\hbox to0pt{\hskip#1 \Einheit
          \raise-4pt\hbox to0pt{\hss\Large$\circ$\hss}\hss}}
\def\Diagonale(#1,#2)#3{\unskip\leavevmode
  \xcoord#1\relax \ycoord#2\relax
      \raise\ycoord \Einheit\hbox to0pt{\hskip\xcoord \Einheit
         \unitlength\Einheit
         \line(1,1){#3}\hss}}
\def\AntiDiagonale(#1,#2)#3{\unskip\leavevmode
  \xcoord#1\relax \ycoord#2\relax \advance\xcoord by -0.05\relax
      \raise\ycoord \Einheit\hbox to0pt{\hskip\xcoord \Einheit
         \unitlength\Einheit
         \line(1,-1){#3}\hss}}
\def\Pfad(#1,#2),#3\endPfad{\unskip\leavevmode
  \xcoord#1 \ycoord#2 \thicklines\ZeichnePfad#3\endPfad\thinlines}
\def\ZeichnePfad#1{\ifx#1\endPfad\let\next\relax
  \else\let\next\ZeichnePfad
    \ifnum#1=1
      \raise\ycoord \Einheit\hbox to0pt{\hskip\xcoord \Einheit
         \vrule height\Pfadd@cke width1 \Einheit depth\Pfadd@cke\hss}%
      \advance\xcoord by 1
     \else\ifnum#1=2
      \raise\ycoord \Einheit\hbox to0pt{\hskip\xcoord \Einheit
         \unitlength\Einheit
         \line(0,1){1}\hss}
      \advance\xcoord by 0
      \advance\ycoord by 1
 \else\ifnum#1=3
      \raise\ycoord \Einheit\hbox to0pt{\hskip\xcoord \Einheit
         \unitlength\Einheit
         \line(1,1){1}\hss}
      \advance\xcoord by 1
      \advance\ycoord by 1
    \else\ifnum#1=4
      \raise\ycoord \Einheit\hbox to0pt{\hskip\xcoord \Einheit
         \unitlength\Einheit
         \line(1,-1){1}\hss}
      \advance\xcoord by 1
      \advance\ycoord by -1
   \else\ifnum#1=5
      \raise\ycoord \Einheit\hbox to0pt{\hskip\xcoord \Einheit
         \unitlength\Einheit
         \line(2,1){2}\hss}
      \advance\xcoord by 2
      \advance\ycoord by 1
	  \else\ifnum#1=6
      \raise\ycoord \Einheit\hbox to0pt{\hskip\xcoord \Einheit
         \unitlength\Einheit
         \line(2,-1){2}\hss}
      \advance\xcoord by 2
      \advance\ycoord by -1
	  \else\ifnum#1=7
      \raise\ycoord \Einheit\hbox to0pt{\hskip\xcoord \Einheit
         \unitlength\Einheit
         \line(3,1){3}\hss}
      \advance\xcoord by 3
      \advance\ycoord by 1
	  \else\ifnum#1=8
      \raise\ycoord \Einheit\hbox to0pt{\hskip\xcoord \Einheit
         \unitlength\Einheit
         \line(3,-1){3}\hss}
      \advance\xcoord by 3
      \advance\ycoord by -1
    \fi\fi\fi\fi\fi\fi\fi\fi
  \fi\next}
\def\hSSchritt{\leavevmode\raise-.4pt\hbox 
to0pt{\hss.\hss}\hskip.2\Einheit
  \raise-.4pt\hbox to0pt{\hss.\hss}\hskip.2\Einheit
  \raise-.4pt\hbox to0pt{\hss.\hss}\hskip.2\Einheit
  \raise-.4pt\hbox to0pt{\hss.\hss}\hskip.2\Einheit
  \raise-.4pt\hbox to0pt{\hss.\hss}\hskip.2\Einheit}
\def\vSSchritt{\vbox{\baselineskip.2\Einheit\lineskiplimit0pt
\hbox{.}\hbox{.}\hbox{.}\hbox{.}\hbox{.}}}
\def\DSSchritt{\leavevmode\raise-.4pt\hbox to0pt{%
  \hbox to0pt{\hss.\hss}\hskip.2\Einheit
  \raise.2\Einheit\hbox to0pt{\hss.\hss}\hskip.2\Einheit
  \raise.4\Einheit\hbox to0pt{\hss.\hss}\hskip.2\Einheit
  \raise.6\Einheit\hbox to0pt{\hss.\hss}\hskip.2\Einheit
  \raise.8\Einheit\hbox to0pt{\hss.\hss}\hss}}
\def\dSSchritt{\leavevmode\raise-.4pt\hbox to0pt{%
  \hbox to0pt{\hss.\hss}\hskip.2\Einheit
  \raise-.2\Einheit\hbox to0pt{\hss.\hss}\hskip.2\Einheit
  \raise-.4\Einheit\hbox to0pt{\hss.\hss}\hskip.2\Einheit
  \raise-.6\Einheit\hbox to0pt{\hss.\hss}\hskip.2\Einheit
  \raise-.8\Einheit\hbox to0pt{\hss.\hss}\hss}}
\def\SPfad(#1,#2),#3\endSPfad{\unskip\leavevmode
  \xcoord#1 \ycoord#2 \ZeichneSPfad#3\endSPfad}
\def\ZeichneSPfad#1{\ifx#1\endSPfad\let\next\relax
  \else\let\next\ZeichneSPfad
    \ifnum#1=1
      \raise\ycoord \Einheit\hbox to0pt{\hskip\xcoord \Einheit
         \hSSchritt\hss}%
      \advance\xcoord by 1
    \else\ifnum#1=2
      \raise\ycoord \Einheit\hbox to0pt{\hskip\xcoord \Einheit
        \hbox{\hskip-2pt \vSSchritt}\hss}%
      \advance\ycoord by 1
    \else\ifnum#1=3
      \raise\ycoord \Einheit\hbox to0pt{\hskip\xcoord \Einheit
         \DSSchritt\hss}
      \advance\xcoord by 1
      \advance\ycoord by 1
    \else\ifnum#1=4
      \raise\ycoord \Einheit\hbox to0pt{\hskip\xcoord \Einheit
         \dSSchritt\hss}
      \advance\xcoord by 1
      \advance\ycoord by -1
    \fi\fi\fi\fi
  \fi\next}
\def\Koordinatenachsen(#1,#2){\unskip
 \hbox to0pt{\hskip-.5pt\vrule height#2 \Einheit width.5pt depth1 
\Einheit}%
 \hbox to0pt{\hskip-1 \Einheit \xcoord#1 \advance\xcoord by1
    \vrule height0.25pt width\xcoord \Einheit depth0.25pt\hss}}
\def\Koordinatenachsen(#1,#2)(#3,#4){\unskip
 \hbox to0pt{\hskip-.5pt \ycoord-#4 \advance\ycoord by1
    \vrule height#2 \Einheit width.5pt depth\ycoord \Einheit}%
 \hbox to0pt{\hskip-1 \Einheit \hskip#3\Einheit 
    \xcoord#1 \advance\xcoord by1 \advance\xcoord by-#3 
    \vrule height0.25pt width\xcoord \Einheit depth0.25pt\hss}}
\def\Gitter(#1,#2){\unskip \xcoord0 \ycoord0 \leavevmode
  \LOOP\ifnum\ycoord<#2
    \loop\ifnum\xcoord<#1
      \raise\ycoord \Einheit\hbox to0pt{\hskip\xcoord 
\Einheit\Punkt\hss}%
      \advance\xcoord by1
    \repeat
    \xcoord0
    \advance\ycoord by1
  \REPEAT}
\def\Gitter(#1,#2)(#3,#4){\unskip \xcoord#3 \ycoord#4 \leavevmode
  \LOOP\ifnum\ycoord<#2
    \loop\ifnum\xcoord<#1
      \raise\ycoord \Einheit\hbox to0pt{\hskip\xcoord 
\Einheit\Punkt\hss}%
      \advance\xcoord by1
    \repeat
    \xcoord#3
    \advance\ycoord by1
  \REPEAT}
\def\Label#1#2(#3,#4){\unskip \xdim#3 \Einheit \ydim#4 \Einheit
  \def\lo{\advance\xdim by-.5 \Einheit \advance\ydim by.5 \Einheit}%
  \def\llo{\advance\xdim by-.25cm \advance\ydim by.5 \Einheit}%
  \def\loo{\advance\xdim by-.5 \Einheit \advance\ydim by.25cm}%
  \def\o{\advance\ydim by.25cm}%
  \def\ro{\advance\xdim by.5 \Einheit \advance\ydim by.5 \Einheit}%
  \def\rro{\advance\xdim by.25cm \advance\ydim by.5 \Einheit}%
  \def\roo{\advance\xdim by.5 \Einheit \advance\ydim by.25cm}%
  \def\l{\advance\xdim by-.30cm}%
  \def\r{\advance\xdim by.30cm}%
  \def\lu{\advance\xdim by-.5 \Einheit \advance\ydim by-.6 \Einheit}%
  \def\llu{\advance\xdim by-.25cm \advance\ydim by-.6 \Einheit}%
  \def\luu{\advance\xdim by-.5 \Einheit \advance\ydim by-.30cm}%
  \def\u{\advance\ydim by-.30cm}%
  \def\ru{\advance\xdim by.5 \Einheit \advance\ydim by-.6 \Einheit}%
  \def\rru{\advance\xdim by.25cm \advance\ydim by-.6 \Einheit}%
  \def\ruu{\advance\xdim by.5 \Einheit \advance\ydim by-.30cm}%
  #1\raise\ydim\hbox to0pt{\hskip\xdim
     \vbox to0pt{\vss\hbox to0pt{\hss$#2$\hss}\vss}\hss}%
}
\catcode`\@=12

\begin{document}
\newtheorem{theorem}{Theorem}
\newtheorem{defn}[theorem]{Definition}
\newtheorem{lemma}[theorem]{Lemma}
\newtheorem{prop}[theorem]{Proposition}
\newtheorem{cor}[theorem]{Corollary}

\begin{center}
{\Large
Card deals, lattice paths, abelian words \\ and combinatorial identities       \\ 
}

\vspace{12mm}
DAVID CALLAN  \\
{\small Dept. of Statistics,  
University of Wisconsin-Madison,
1300 University Ave,
Madison, WI \ 53706}  \\

\vspace*{2mm}

{\bf callan@stat.wisc.edu}  \\
\vspace{2mm}

December 26, 2008   
\end{center}

\begin{abstract}
We give combinatorial interpretations of 
several related identities associated with the names Barrucand, Strehl 
and Franel, including one for the  Ap\'{e}ry numbers,  
$\sum_{k=0}^{n}\binom{n}{k}\binom{n+k}{k}\sum_{j=0}^{k}\binom{k}{j}^{3}=
\sum_{k=0}^{n}\binom{n}{k}^{2}\binom{n+k}{k}^{2}$. The combinatorial 
constructs employed are derangement-type card 
deals as introduced in a previous paper on Barrucand's identity, labeled 
lattice paths and, following a comment of Jeffrey Shallit, abelian words over a 3-letter alphabet.
\end{abstract}

\section{Introduction} The purpose of this paper is to
give simple direct combinatorial interpretations of two identities of Strehl \cite{twostrehl}, for 
the 
\htmladdnormallink{Franel}{http://www.research.att.com:80/cgi-bin/access.cgi/as/njas/sequences/eisA.cgi?Anum=A000172}
and  
\htmladdnormallink{Ap\'{e}ry}{http://www.research.att.com:80/cgi-bin/access.cgi/as/njas/sequences/eisA.cgi?Anum=A005259}
numbers respectively,
\begin{equation}
    \sum_{k=0}^{n}\binom{n}{k}^{3}=
\sum_{k=0}^{n}\binom{n}{k}^{2}\binom{2k}{n},
    \label{franel}
\end{equation}
and
\begin{equation}
\sum_{k=0}^{n}\binom{n}{k}\binom{n+k}{k}\sum_{j=0}^{k}\binom{k}{j}^{3}  = 
    \sum_{k=0}^{n}\binom{n}{k}^{2}\binom{n+k}{k}^{2},
    \label{apery}
\end{equation}
and of the following curious sequence of identities involving powers of 
successively larger integers,
\begin{equation}
    \sum_{k=0}^{n} \binom{n}{k}\binom{2k}{k}2^{k} =
\sum_{k=0}^{n} \binom{n}{k}\binom{2n-k}{n}3^{k} =
\sum_{k=0}^{n}\binom{n}{k}^{2}4^{k}  =
\sum_{k=0}^{n/2} \binom{n}{2k}\binom{2k}{k}4^{k}5^{n-2k}.
    \label{34}
\end{equation}
The first three of these expressions are equated in \cite[Eqs. 
34,\:35]{strehl94}, and all give sequence
\htmladdnormallink{A084771}{http://www.research.att.com:80/cgi-bin/access.cgi/as/njas/sequences/eisA.cgi?Anum=A084771}
in OEIS.

The combinatorial constructs employed are (generalizations of) the derangement-type card 
deals introduced in a previous paper on Barrucand's identity \cite{barr}, 
the labeled lattice paths cited by Nour-Eddine Fahssi in
\htmladdnormallink{A084771}{http://www.research.att.com:80/cgi-bin/access.cgi/as/njas/sequences/eisA.cgi?Anum=A084771},
and, following a comment of Jeffrey Shallit \cite{shallitcomm}, 
abelian words over a 3-letter alphabet.

Section 2 reviews the card deals and abelian words/matrices. Section 3 
presents a 1-to-1 correspondence between them
and reinterprets Barrucand's identity,
\begin{equation}
    \sum_{k=0}^{n}\binom{n}{k}\sum_{j=0}^{k}\binom{k}{j}^{3}=
\sum_{k=0}^{n}\binom{n}{k}^{2}\binom{2k}{k} 
    \label{barrucand}
\end{equation}
in terms of abelian matrices. Section 4 gives interpretations for 
(\ref{franel}) and Section 5 for (\ref{apery}). Section 6 presents 
three equinumerous combinatorial constructs involving lattice paths, card deals and 
matrices  respectively, and Section 7 uses 
them to interpret (\ref{34}). 

\section{Card deals and abelian words/matrices} \label{deals}
A \emph{Barrucand $n$-deal} \cite{barr} is formed as follows. 
Start with a deck of $3n$ cards, $n$ each colored red, green and 
blue, in denominations 1 through $n$,
choose an arbitrary subset of the denominations
and  deal all cards of the chosen denominations 
into three equal-size hands to players designated red, green and 
blue in such a way that no player receives a card of her own color.
Let $\bb_{n}$ denote the set of Barrucand $n$-deals.

The left side of (\ref{barrucand}) 
counts $\bb_{n}$ by total number of cards, $k$, in red's hand and number of green 
cards, $j$, in red's hand: first, there are $\binom{n}{k}$ ways to choose the 
denominations in the deal; next, $j$ green cards in red's hand 
implies both 
$j$ blue cards in green's hand and $j$ red cards in blue's hand, and 
these cards determine the deal. For each hand there are $\binom{k}{j}$ 
ways to choose the determining cards, so $\binom{k}{j}^{3}$ choices in all. 
As shown in \cite{barr}, the right 
side counts $\bb_{n}$ by number of distinct denominations, $k$, in red's 
hand; another approach to establishing this count is given below.

Serendipitously, on the day \cite{barr} was published, the editor 
emailed me that the counting sequence for $\bb_{n}$ also arose in his 
recently posted paper \cite{abelian} counting abelian squares. An \emph{abelian 
square} (over an alphabet) is a word of the form $ww'$ where $w'$ is a 
rearrangement of $w$. Its \emph{size} is the number of letters in $w$ (= 
number of letters in $w'$). As easily seen, the number of abelian squares over a 
three-letter alphabet, say $\{1,2,3\}$, of size $n$ with 
$n-k$ 1s in $w$ is $\binom{n}{k}^{2}\binom{2k}{k}$ \cite{abelian}, 
the summand on the right in (\ref{barrucand}). This raises the 
questions of a bijection from $\bb_{n}$ to abelian squares over $\{1,2,3\}$ 
and of an abelian squares interpretation for the left side of (\ref{barrucand}). 
It is convenient to represent an abelian square $ww'$ of size $n$  as a 
$2 \times n$
matrix $\left(
\begin{smallmatrix}
    w \\
    w'
\end{smallmatrix} \right)$, a so-called \emph{abelian matrix}, 
so that we can refer to its columns.

\section{Bijection from Barrucand deals to abelian matrices} 
\label{hope}
The following table describes a bijection from $\bb_{n}$, the set of 
Barrucand $n$-deals, to $2\times n$ abelian 
matrices over $\{1,2,3\}$ by specifying the locations of the 9 
possible distinct columns in the matrix ($R,G,B$ are short for red, green, blue respectively).
\begin{center}
\begin{tabular}{|c|c|}
    \hline
  \raisebox{-0.6ex}[0pt]{matrix}  & \\
  \raisebox{.2ex}[0pt]{column} & \raisebox{1.3ex}[0pt]{locations given 
  by denominations that are \ldots} \\ \hline\hline
   & \\[-10pt]
  {\normalsize $\begin{smallmatrix}
        1 \\ 1
    \end{smallmatrix}$} &  in $[n]$, not in deal \\[5pt]
  $\begin{smallmatrix}
      1 \\ 2
  \end{smallmatrix}$ &  in deal, not in red's hand and not on R 
  in blue's hand \\[5pt]
   $\begin{smallmatrix}
      1 \\ 3
  \end{smallmatrix}$ &  not in red's hand but do occur on R 
  in blue's hand \\[5pt] \hline
   & \\[-10pt]
   $\begin{smallmatrix}
      2 \\ 1
  \end{smallmatrix}$ &  in red's hand on G and B and also occur on R 
  in blue's hand \\[5pt]
   $\begin{smallmatrix}
      2 \\ 2
  \end{smallmatrix}$ &  in red's hand on G only and also occur on R 
  in blue's hand \\[5pt]
   $\begin{smallmatrix}
      2 \\ 3
  \end{smallmatrix}$ &  in red's hand on B only and also occur on R 
  in blue's hand \\[5pt] \hline
   & \\[-10pt]
  $\begin{smallmatrix}
      3 \\ 1
  \end{smallmatrix}$ &  in red's hand on G and B and don't occur on R 
  in blue's hand \\[5pt]
   $\begin{smallmatrix}
      3 \\ 2
  \end{smallmatrix}$ &  in red's hand on G only and don't occur on R 
  in blue's hand \\[5pt]
   $\begin{smallmatrix}
      3 \\ 3
  \end{smallmatrix}$ &  in red's hand on B only and don't occur on R 
  in blue's hand \\[5pt]  \hline
\end{tabular}\\[4mm]
Bijection from deals to matrices \\[1mm]
Table 1
\end{center}
Note, for example, that the denominations not in red's hand give the locations of 1s in the top row. 
It is straightforward to check that this mapping is a bijection as 
claimed and that its inverse is given by the following table.
\begin{center}
\begin{tabular}{|c|cc|} 
    \hline
     &  \raisebox{-0.6ex}[0pt]{denominations}  &  \raisebox{-0.6ex}[0pt]{given by} \\     
\raisebox{1.3ex}[0pt]{ \quad player\quad} & \raisebox{.2ex}[0pt]{on \ldots 
cards} & \raisebox{.2ex}[0pt]{\quad locations 
  of \ldots \quad} \\ \hline\hline
   & & \\[-10pt]
  & G and B & $\begin{smallmatrix}
        2 \\ 1
    \end{smallmatrix}$, $\begin{smallmatrix}
        3 \\ 1
    \end{smallmatrix}$ \\[5pt]
  red & G only &$\begin{smallmatrix}
        2 \\ 2
    \end{smallmatrix}$,  $\begin{smallmatrix}
        3 \\ 2
    \end{smallmatrix}$ \\[5pt]
  & B only &  $\begin{smallmatrix}
        2 \\ 3
    \end{smallmatrix}$,  $\begin{smallmatrix}
        3 \\ 3
    \end{smallmatrix}$ \\[5pt] \hline
  & & \\[-10pt]    
    & B and R & $\begin{smallmatrix}
        1 \\ 2
    \end{smallmatrix}$, $\begin{smallmatrix}
        3 \\ 2
    \end{smallmatrix}$ \\[5pt]
  green & B only &$\begin{smallmatrix}
        1 \\ 3
    \end{smallmatrix}$,  $\begin{smallmatrix}
        2 \\ 2
    \end{smallmatrix}$ \\[5pt]
  & R only &  $\begin{smallmatrix}
        3 \\ 1
    \end{smallmatrix}$,  $\begin{smallmatrix}
        3 \\ 3
    \end{smallmatrix}$ \\[5pt] \hline
    & & \\[-10pt]     
    & R and G & $\begin{smallmatrix}
        1 \\ 3
    \end{smallmatrix}$, $\begin{smallmatrix}
        2 \\ 3
    \end{smallmatrix}$ \\[5pt]
  blue & R only &$\begin{smallmatrix}
        2 \\ 1
    \end{smallmatrix}$,  $\begin{smallmatrix}
        2 \\ 2
    \end{smallmatrix}$ \\[5pt]
  & G only &  $\begin{smallmatrix}
        1 \\ 2
    \end{smallmatrix}$,  $\begin{smallmatrix}
        3 \\ 3
    \end{smallmatrix}$ \\[5pt] \hline

\end{tabular}\\[4mm]
Bijection from matrices to deals \\[1mm]
Table 2
\end{center}
For example, with $n=5$ and subscripts referring to card color, the deal for which red's hand contains 
$2_{G},\,2_{B},\,4_{B},\,5_{G}$, green's hand  contains 
$1_{B},\,2_{R},\,4_{R},\,5_{B}$, and blue's hand contains
$1_{G},\,1_{R},\,4_{G},\,5_{R}$ corresponds to the abelian matrix
$\left(
\begin{smallmatrix}
   1 & 3 & 1 & 3 & 2 \\
   3 & 1 & 1 & 3 & 2
\end{smallmatrix}
\right) $. 

Evidently, abelian matrices are somewhat more concise than 
Barrucand deals but, on the other hand, some statistics on 
$\bb_{n}$ are more appealing than their counterparts for abelian matrices.
For example,

\vspace*{-3mm}

\begin{center}
    \begin{tabular}{ccc}
        
        \#\,cards in red's hand & $\leftrightarrow$ & $n- \# \,\left( \begin{smallmatrix}
   1  \\
   1
\end{smallmatrix}\right)$ 
columns  \\
        
        \#\,distinct denominations in red's hand & $\leftrightarrow$ & 
        total \#\ 2s and 3s in top row  \\
       
        \#\,green cards in red's hand & $\leftrightarrow$ & \# columns
$\left( \begin{smallmatrix}
   p  \\
   q
\end{smallmatrix}\right)$ 
with $p>1$ and $q<3$.  \\
        
    \end{tabular}
\end{center}

\vspace*{-2mm}

In particular, using these correspondences and the second paragraph 
of Section \ref{deals}, the left side of Barrucand's identity (\ref{barrucand}) 
counts abelian matrices of size $n$ over $\{1,2,3\}$ by number, $k$, of 
columns
$\left( \begin{smallmatrix}
   p  \\
   q
\end{smallmatrix}\right) \ne 
\left( \begin{smallmatrix}
   1  \\
   1
\end{smallmatrix}\right)$   
and number, $j$, of columns
$\left( \begin{smallmatrix}
   p  \\
   q
\end{smallmatrix}\right)$ 
with $p>1$ and $q<3$. Summarizing these observations, we have the 
following alternative interpretation.
\begin{prop}\label{abelianword}
    For Barrucand's identity (\ref{barrucand}), the right side of counts abelian words $ww'$ of 
    length $2n$ by number, $n-k$, of 1s in $w$ while the left side 
    counts them by number of positions, $n-k$, in which both $w$ and 
    $w'$ have a 1.
\end{prop}

A generalization of Barrucand's identity (identity (37) in 
\cite{strehl94}),
\begin{equation} \label{genbarr}
 \sum_{k=0}^{n}\binom{n}{k}\sum_{j=0}^{k}\binom{k}{j}^{2}\binom{k}{j-a}=
\sum_{k=0}^{n}\binom{n}{k}^{2}\binom{2k}{k-a}, 
\end{equation}
can be treated similarly. Let $\a_{n,a}$ denote the set of $2\times 
n$ matrices with entries in $\{1,2,3\}$, the same number of 1s in each 
row, and $a$ more 3s in the top row than in the bottom row. For 
example, 
$\left( \begin{smallmatrix}
   1 & 2 & 3  \\
   2 & 1 & 2
\end{smallmatrix}\right) \in \a_{3,1}$, and $a=0$ gives abelian 
matrices. 
Then the two sides of (\ref{genbarr}) count $\a_{n,a}$ by the very 
same statistics as the two sides of (\ref{barrucand}) count abelian 
matrices.

\section[Franel numbers]{Franel numbers, $\protect\mbf{\sum_{k=0}^{n}\binom{n}{k}^{3}=
\sum_{k=0}^{n}\binom{n}{k}^{2}\binom{2k}{n}}$}
A \emph{Franel $n$-deal} is a Barrucand $n$-deal in which \emph{all} the cards are 
dealt to the players. Let $\f_{n}$ denote the set of Franel $n$-deals. 
As observed in Section 2, the left side of the identity for the 
\htmladdnormallink{Franel}{http://www.research.att.com:80/cgi-bin/access.cgi/as/njas/sequences/eisA.cgi?Anum=A000172}
numbers counts $\f_{n}$ by number,
$k$, of green cards in red's hand. Translated to abelian matrices, the left side 
counts $\f_{n}'$, the abelian matrices of size $n$ over $\{1,2,3\}$ with no 
$\left( \begin{smallmatrix}
   1  \\
   1
\end{smallmatrix}\right)$ 
columns, by number, $k$, 
of columns
$\left( \begin{smallmatrix}
   p  \\
   q
\end{smallmatrix}\right)$ 
with $p>1$ and $q<3$.

As for the right side, let us count $\f_{n}'$ by number, 
$j$, of 1s in each row: $\binom{n}{j}$ [place 1s in top row] $\times 
\binom{n-j}{j}$ [place 1s in bottom row] $\times \binom{2n-2j}{n-j}$ 
[choose $n-j$ of the remaining $2n-2j$ positions; place 2s in the chosen positions 
in the top row and fill out the top row with 3s; place 3s in the chosen positions 
in the bottom row and fill out the bottom row with 2s]. (The latter clever 
argument is due to Richmond and Shallit \cite{abelian}.) 
Thus, with $k:=n-j$, the number  of abelian matrices in $\f_{n}'$
with a total of $k$ 2s and 3s in each row is 
$\binom{n}{n-k}\binom{k}{n-k}\binom{2k}{k}=\binom{n}{k}^{2}\binom{2k}{n}$. 
Translated back to card deals, the right side counts $\f_{n}$ by 
number of distinct denominations in red's hand.

\section[Apery numbers]{Ap\'{e}ry numbers, \\ $\protect\mbf{\sum_{k=0}^{n}\binom{n}{k}\binom{n+k}{k}\sum_{j=0}^{k}\binom{k}{j}^{3}  = 
    \sum_{k=0}^{n}\binom{n}{k}^{2}\binom{n+k}{k}^{2}}$}
The counting sequence for this identity,
\htmladdnormallink{A005259}{http://www.research.att.com:80/cgi-bin/access.cgi/as/njas/sequences/eisA.cgi?Anum=A005259},
cropped up in Roger Ap\'{e}ry's celebrated proof of the irrationality 
of $\zeta(3)$ \cite{poorten78} and the identity inspired a survey paper by 
Volker Strehl \cite{strehl94} in which he offers six different proofs 
including a combinatorial proof of a substantial generalization
and, indeed, proves most of the other identities in this paper.
Still, simple direct fully bijective proofs may be of interest.

Let $\bb_{n,k}$ denote the set of deals in $\bb_{n}$ with $k$ cards in 
red's hand, equivalently, $k$ denominations 
in the deal. Thus $\v\,\bb_{n,k}\,\v=\binom{n}{k}\sum_{j=0}^{k}\binom{k}{j}^{3}$. 
To get the left side of Ap\'{e}ry (\ref{apery}), 
we need an additional factor of $\binom{n+k}{k}$ on the left side of 
Barrucand (\ref{barrucand}). This motivates us to consider a simple 
construction and define $\a_{n,k}$ to be the set of 
pairs $(D,i)$ where $D \in \bb_{n,k}$ and $1 \le i \le \binom{n+k}{k}$. 
Thus $\v\,\a_{n,k}\,\v=\binom{n+k}{k} \v\,\bb_{n,k}\,\v= 
\binom{n}{k}\binom{n+k}{k}\sum_{j=0}^{k}\binom{k}{j}^{3}$ and 
$\a_{n}:=\bigcup_{k=0}^{n}\a_{n,k}$ is counted by the left side of Ap\'{e}ry.
\begin{prop}
    Just as for Barrucand, the right side of Ap\'{e}ry counts 
    $\a_{n}$ by number of distinct denominations in the red player's 
    hand in the associated deal.
    \label{prop1}
\end{prop}
The proof needs the identity
\begin{equation}
    \sum_{a\ge 0}\binom{k}{a}\binom{n-k}{a}\binom{n+k+a}{n} = 
    \binom{n+k}{k}\binom{n+k}{n-k},
    \label{andrewsid}
\end{equation}
proved combinatorially by George Andrews \cite{andrews} in a more 
general form (see also \cite[Eqs. (19) and (20)]{strehl94}). Applied to 
(\ref{andrewsid}), his proof shows that the right side counts pairs 
$(K,L)$ where $K$ is a $k$-element subset of $[n+k]$ and $L$ is an 
$(n-k)$-element subset of  $[n+k]$ while the left side counts these 
pairs by ``intermingling coefficient'' $a$: the number of elements 
in $L$ among the $k$ smallest elements of $K\cup L$.

A proof of Prop. \ref{prop1} can now be devised following the 
analysis of $\bb_{n}$ in \cite{barr} but it is a little simpler to 
translate to abelian matrices and prove the following equivalent 
result.  
\begin{prop}
Let $\a_{n}'$ denote the set of pairs $(A,i)$ with $A$ a 
$2\times n$ abelian matrix over $\{1,2,3\}$ and $1 \le i \le 
\binom{n+j}{j}$ where $n-j$ is the number of 
$\left( \begin{smallmatrix}
   1  \\
   1
\end{smallmatrix}\right)$ 
columns in $A$.

Then $\sum_{k=0}^{n}\binom{n}{k}^{2}\binom{n+k}{k}^{2}$ counts 
$\a_{n}'$ by total number, $k$, of $2$s and $3$s in the top row.
\end{prop}
\noindent \textbf{Proof} \quad 
Suppose $(A,i)\in\a_{n}'$ has $k$ 2s and 3s, hence $n-k$ 1s, in the 
top row. Now count by number of 
$\left( \begin{smallmatrix}
   1  \\
   1
\end{smallmatrix}\right)$ 
columns, say $n-k-a$. Thus we have $\binom{n+k+a}{k+a}$ choices for 
the second member $i$ of the pair $(A,i)$ and choices for $A$ as 
follows: place 1s in top row [\,$\binom{n}{n-k}$ choices\,],
locate 
$\left( \begin{smallmatrix}
   1  \\
   1
\end{smallmatrix}\right)$
columns [\,$\binom{n-k}{n-k-a}$ choices\,], place $a$ 1s in the bottom row not 
below 1s in the top row [\,$\binom{k}{a}$ choices\,], place 2s 
and 3s [\,$\binom{2k}{k}$ choices, as explained in Section 4\,].
All told, the number of choices for $(A,i)$ is 
\begin{eqnarray*}
\binom{n}{k}\binom{2k}{k}\sum_{a\ge 
0}\binom{n-k}{a}\binom{k}{a}\binom{n+k+a}{n} & = &  
     \binom{n}{k}\binom{2k}{k}\binom{n+k}{k}\binom{n+k}{n-k}  \\
     & = & \binom{n}{k}^{2}\binom{n+k}{k}^{2},
\end{eqnarray*}
using (\ref{andrewsid}) at the first equality. \qed

\section[Combinatorial constructs]{Combinatorial constructs for  (\ref{34})} 
A \emph{Delannoy path} is a lattice path of upsteps $U=(1,1)$, downsteps 
$D=(1,-1)$, and flatsteps $F=(1,0)$ with an equal number of $U$s and $D$s. 
The line joining its endpoints, necessarily horizontal, is \emph{\gl}. Each upstep in a 
Delannoy path has a matching downstep (and conversely): given an 
upstep above \gl (resp. below \gl), travel directly east (resp. west) 
until you encounter a downstep.

\vspace*{-5mm}

\Einheit=0.8cm
\[
\Pfad(-7,2),1314441343334\endPfad
\SPfad(-7,2),111111111111111\endSPfad
\Label\o{ \overleftrightarrow{\phantom{aaaaa}}}(-4.5,2.1)
\Label\o{ \overleftrightarrow{\phantom{aaaaaaaaaaaaaaaaaaaa}}}(0.5,1.1)
\Label\o{ \overleftrightarrow{\phantom{aaaaaa}}}(-0.5,0.2)
\Label\o{ \overleftrightarrow{\phantom{aa}}}(2,0.2)
\Label\o{ \overleftrightarrow{\phantom{aa}}}(5,2.0)
\DuennPunkt(-7,2)
\DuennPunkt(-6,2)
\DuennPunkt(-5,3)
\DuennPunkt(-4,3)
\DuennPunkt(-3,2)
\DuennPunkt(-2,1)
\DuennPunkt(-1,0)
\DuennPunkt(0,0)
\DuennPunkt(1,1)
\DuennPunkt(2,0)
\DuennPunkt(3,1)
\DuennPunkt(4,2)
\DuennPunkt(5,3)
\DuennPunkt(6,2)
\Label\o{\textrm{{\footnotesize ground}}}(7,2.0)
\Label\u{\textrm{{\footnotesize level}}}(7,2.1)
\Label\u{ \textrm{matching step pairs in a Delannoy path}}(0,-.6)
\]

\vspace*{2mm}

\noindent Thus the slanted steps ($U$ and $D$) in a Delannoy path are partitioned 
into matching pairs of opposite-slope steps.

A \emph{Hanna $n$-path}
is a Delannoy path with $n$ labeled steps: each slanted step gets one 
of two labels (colors), say 1 or 2, and each flat step gets one 
of five labels, say 1,\,2,\,3,\,4 or 5. As observed by Nour-Eddine 
Fahssi, Hanna $n$-paths are counted by 
\htmladdnormallink{A084771}{http://www.research.att.com:80/cgi-bin/access.cgi/as/njas/sequences/eisA.cgi?Anum=A084771}.
 
A \emph{Hanna $n$-deal} is formed in the same way as a Barrucand deal 
except that the 
hands need not all be of equal size: if there are $j$ denominations 
in the deal, only red's hand is required to contain its fair share of $j$ cards and the 
remaining $2j$ cards are split arbitrarily between the green 
and blue players.

A \emph{Hanna $n$-matrix} is a $2\times n$ matrix with entries in $\{1,2,3\}$ and the 
same number of 1s in each row.

Hanna $n$-matrices, $n$-deals, and $n$-paths are equinumerous:
the mapping in Table 1 of Section \ref{hope} (with a larger domain) is a bijection from the matrices 
to the deals, and there is a simple bijection from the matrices to 
the paths: transform each column in turn (subscripts denote step 
labels) according to the following table.

\vspace*{-6mm}

\[
\begin{array}{c|ccccccccc}
    \textrm{matrix column} &
    \begin{smallmatrix}
   1  \\
   1
\end{smallmatrix} & \begin{smallmatrix}
   1  \\
   2
\end{smallmatrix} & \begin{smallmatrix}
   1  \\
   3
\end{smallmatrix} & \begin{smallmatrix}
   2  \\
   1
\end{smallmatrix} & \begin{smallmatrix}
   2  \\
   2
\end{smallmatrix} & \begin{smallmatrix}
   2  \\
   3
\end{smallmatrix} & \begin{smallmatrix}
   3  \\
   1
\end{smallmatrix} & \begin{smallmatrix}
   3  \\
   2
\end{smallmatrix} & \begin{smallmatrix}
   3  \\
   3
\end{smallmatrix}   \\[2mm] \hline 
\rule[-0.4cm]{0mm}{1cm}\textrm{labeled step} &  F_{1} & U_{1} & U_{2} & D_{1} & F_{2} & F_{3} & 
D_{2} & F_{4} & F_{5} 
\end{array}
\]

\vspace*{-3mm}

In the next section, we use these constructs to give a combinatorial 
interpretation of the identities (\ref{34}).

\section[Combinatorial interpretations]{Combinatorial interpretations for  (\ref{34})}

The summand in the first expression in (\ref{34}), 
$\binom{n}{k}\binom{2k}{k}2^{k}$, is the number of Hanna 
$n$-deals with $k$ cards in red's hand. To see this, expand 
$2^{k}$ as $\sum_{j=0}^{k}\binom{k}{j}$. Then the resulting summand,
$\binom{n}{k}\binom{k}{j}\binom{2k}{k}$, is the number of  Hanna 
$n$-deals with $k$ cards in red's hand and $j$ red cards in blue's hand: 
choose denominations in the deal [\,$\binom{n}{k}$ choices\,], choose red 
denominations in blue's hand [\,$\binom{k}{j}$ choices\,] and the remaining red 
cards are forced into green's hand, select red's hand from the green 
and blue cards [\,$\binom{2k}{k}$ choices\,] and the remaining green and blue 
cards are forced into the hand of opposite color. \qed

The least obvious statistic for the sums in (\ref{34}) is the 
one for the second sum. Actually, it is a sum of two statistics. On 
Hanna $n$-paths, define the statistic $X$ to be the number of matching 
pairs of slanted steps not both labeled 1, and define $Y$ to be the 
number of flatsteps whose label exceeds 2. Then the summand in the 
second expression in (\ref{34}), $\binom{n}{k}\binom{2n-k}{n}3^{k}$, is the 
number of Hanna $n$-paths for which $X+Y=k$. This is an immediate 
consequence of the following two propositions.
\begin{prop} \label{a}
    The number of Hanna $n$-paths with $X=i$ and $Y=j$ is 
    \[
    \binom{n}{j}\binom{n-j}{i}\binom{2n-2i-2j}{n-j}3^{i+j}.
    \]
\end{prop}
\begin{prop} \label{b}
    \[
  \sum_{\substack{i,\,j: \\ i+j=k}} 
  \binom{n}{j}\binom{n-j}{i}\binom{2n-2i-2j}{n-j}3^{i+j}= \binom{n}{k}\binom{2n-k}{n}3^{k}.
\]
\end{prop}

\noindent \textbf{Proof of Prop. (\ref{a})}\quad To form a Hanna $n$-path with $X=i$ and $Y=j$, 
choose locations in the path for flatsteps whose label exceeds 2 
[\,$\binom{n}{j}$ choices\,], label these flatsteps [\,$3^{j}$ choices\,], 
choose locations for the upsteps in matching pairs whose members are 
not both labeled 1 [\,$\binom{n-j}{i}$ choices\,], assign labels to these 
pairs [\,$3^{i}$ choices, since each $U$-$D$ pair may be labeled 1-2, 
2-1, or 2-2]. Now consider the steps in the $n-i-j$ locations not yet 
filled (including the downsteps in the matching pairs). These steps 
form a path of $U$s, $D$s, and $F$s of length $n-i-j$ with $i$ more 
$D$s than $U$s. The labels on the slanted steps in this path are already 
determined and the flatsteps are bicolored (labeled 1 or 2). Expanding 
the path via the transformation rules $U\rightarrow UU,\ D\rightarrow 
DD,\ F_{1}\rightarrow UD,\ F_{2}\rightarrow DU$ (subscript denotes 
label), it becomes a  path of $U$s and $D$s of length $2n-2i-2j$ with 
$n-2i-j\ U$s and $n-j\ D$s. There are $\binom{2n-2i-2j}{n-j}$ such 
paths, and the expansion is reversible. Thus all factors in the expression of Prop. (\ref{a}) have been 
accounted for. \qed

\noindent \textbf{Proof of Prop. (\ref{b})}
\begin{eqnarray*}
\sum_{\substack{i,\,j: \\ i+j=k}} 
  \binom{n}{j}\binom{n-j}{i}\binom{2n-2i-2j}{n-j}3^{i+j} & = &  
     \sum_{j}\binom{n}{j}\binom{n-j}{k-j}\binom{2n-2k}{n-j}3^{k}  \\
     & = & \sum_{j}\binom{n}{k}\binom{k}{j}\binom{2n-2k}{n-j}3^{k} \\
     & = & \binom{n}{k}\binom{2n-k}{n}3^{k},
\end{eqnarray*}
using the Chu-Vandermonde identity at the last equality. \qed

The summand in the third expression in (\ref{34}), $\binom{n}{k}^{2}4^{k}$, is the 
number of $2\times n$ Hanna $n$-matrices with $n-k$ 1s in each row: place 
the 1s [\,$\binom{n}{n-k}^{2}=\binom{n}{k}^{2}$ choices\,] and then fill the remaining $2k$ entries 
with 2s and 3s arbitrarily [\,$2^{2k}$ choices\,]. 
Equivalently, it counts Hanna $n$-deals by number, $n-k$, of denominations appearing in 
red's hand. (Alternative interpretations of the other expressions in 
(\ref{34}) are left to the reader.) \qed

The summand in the fourth expression in (\ref{34}), $\binom{n}{2k}\binom{2k}{k}4^{k}5^{n-2k}$, is the 
number of Hanna $n$-paths with $k$ upsteps: choose locations for the 
slanted steps [\,$\binom{n}{2k}$ choices\,], insert $U$s and $D$s into 
these locations [\,$\binom{2k}{k}$ choices\,], label the slanted steps 
[\,$2^{2k}$ choices\,], and lastly, label the flatsteps [\,$5^{n-2k}$ 
choices\,].

\end{document}